\documentclass{amsart}
\usepackage{helvet}
\usepackage{comment}
\usepackage{setspace}
\usepackage{enumitem}
\usepackage[all, cmtip]{xy}
\usepackage{amssymb}
\usepackage{graphicx}
\usepackage{ragged2e}
\newcommand{\Proj}{\mathrm{Proj}}
\newcommand{\Flat}{\mathrm{Flat}}
\newcommand{\Cot}{\mathrm{Cot}}
\newcommand{\Inj}{\mathrm{Inj}} 
\newcommand{\FPInj}{\mathrm{FP\text{-}injective}}
\newcommand{\GF}{\mathcal{GF}}
\newcommand{\PGF}{\mathcal{PGF}}
\newcommand{\GP}{\mathcal{GP}}
\newcommand{\DP}{\mathcal{DP}}
\newcommand{\DI}{\mathcal{DI}}
\newcommand{\wDI}{\emph{w}\mathcal{DI}}
\newcommand{\GI}{\mathcal{GI}}
\newcommand{\Gpd}{\mathrm{Gpd}}
\newcommand{\Gfd}{\mathrm{Gfd}}
\newcommand{\Gid}{\mathrm{Gid}}
\newcommand{\pd}{\mathrm{pd}}
\newcommand{\id}{\mathrm{id}}
\newcommand{\fd}{\mathrm{fd}}
\newcommand{\DPdim}{\mathrm{DP\text{-}dim}}
\newcommand{\PGFdim}{\mathrm{PGF\text{-}dim}}
\usepackage[
    style=numeric,
    sorting=nyt,  
    maxnames=5,
    minnames=1,
    giveninits=true,
    doi=true,
    url=true,
    backend=biber
]{biblatex}
\DefineBibliographyExtras{english}{

}

\DeclareNameAlias{sortname}{family-given}
\DeclareNameAlias{default}{family-given}

\AtBeginBibliography{
  
}
\DeclareFieldFormat[article]{title}{#1}
\DeclareFieldFormat[book]{title}{\mkbibemph{#1}}
\DeclareFieldFormat{journaltitle}{#1}

\usepackage{xcolor} 
\definecolor{darkblue}{RGB}{0,0,139} 
\usepackage{etoolbox} 
\usepackage[margin=1in]{geometry} 
\usepackage{tikz}
\usepackage{amsmath,amsthm,amssymb}
\usepackage[colorlinks=true,  
            linkcolor=darkblue,  
            citecolor=darkblue, 
            urlcolor=darkblue 
           ]{hyperref}
\usepackage{abstract}

\usepackage{titlesec}
\titleformat{\section}
  {\normalfont\Large\bfseries\color{darkblue}}
  {\thesection}{1em}{}
\titleformat{\subsection}
  {\normalfont\large\bfseries\color{darkblue}}
  {\thesubsection}{1em}{}
\titleformat{\subsubsection}
  {\normalfont\normalsize\bfseries\color{darkblue}}
  {\thesubsubsection}{1em}{}

\newtheoremstyle{plainnormal}  
  {3pt}      
  {3pt}      
  {\normalfont} 
  {}         
  {\color{darkblue}\bfseries} 
  {.}        
  {0.5em}    
  {}         

\theoremstyle{plainnormal}
\newtheorem{theorem}{Theorem}[section]
\newtheorem{lemma}[theorem]{Lemma}
\newtheorem{corollary}[theorem]{Corollary}
\newtheorem{proposition}[theorem]{Proposition}

\newtheorem{remark}[theorem]{Remark}

\makeatletter
\renewenvironment{proof}[1][\proofname]{
  \par\pushQED{\qed}\normalfont
  \topsep6pt \partopsep0pt
  \trivlist
  \item[\hskip\labelsep\color{darkblue}\itshape #1.]
}{
  \popQED\endtrivlist\@endpefalse
}
\makeatother
\usepackage{cleveref}

\begin{document}
\onehalfspacing
\title{On classical and Gorenstein homological invariants of rings}
\author{José Manuel Fresneda}
\date{}

\thanks{J.M.F. was partly supported by grant PID2024-158993NB-100 and by grant PID2024-155576NB-I00. Both grants are funded by MICIU/AEI/10.13039/501100011033 /FEDER, UE}


\maketitle
\begin{center}
\begin{minipage}{0.8\textwidth}
\textbf{\textcolor{darkblue}{Abstract.}}
\justifying We prove that, over any ring $R$, the supremum of the projective dimensions of the flat left $R$-modules coincides with the supremum of the Gorenstein projective dimensions of the Gorenstein flat left $R$-modules. As a consequence, we obtain new characterizations of left $n$-perfect rings in terms of Gorenstein projective, Ding projective, and projectively coresolved Gorenstein flat dimensions, extending results by Emmanouil and Dalezios and by Christensen, Estrada, and Thompson. We also introduce Gorenstein analogues of several classical homological invariants and study their relationships with the classical ones, identifying conditions under which they coincide with the classical invariants. Finally, we obtain characterizations of left weakly $n$-$\Sigma$-cotorsion rings introduced by Cortés-Izurdiaga, Estrada and Fresneda in terms of Gorenstein classes of modules.
\end{minipage}
\end{center}
\section{Introduction}
Homological invariants play a central role in the study of modules over rings. Among the classical invariants associated to a ring $R$ are those defined via homological dimensions of distinguished classes of $R$-modules, such as $\mathrm{splf}(R)$, the supremum of the projective lengths of flat left $R$-modules. Recall that a ring $R$ is said to be \emph{left $n$-perfect} (for some integer $n \geq 0$) if every flat left $R$-module has projective dimension at most $n$, equivalently, if $\mathrm{splf}(R) \leq n$.

The first main goal of this paper is to introduce and study the Gorenstein counterpart of this invariant, denoted by $\mathrm{G\text{-}splf}(R)$, defined as the supremum of the Gorenstein projective lengths of Gorenstein flat left $R$-modules. By employing a generalization of a result of Guil Asensio and Herzog (see \cite[Theorem 19 and Corollary 20]{GuilAsensioHerzog2005}), as established in \cite[Theorem 4.2]{CortesIzurdiagaEstradaFresneda2026}, we prove our main result (see Theorem~\ref{th:s(GP)l(GF)=splf}):

\medskip\par\noindent
\textbf{Theorem A.} For any ring $R$, one has $\mathrm{splf}(R) = \mathrm{G\text{-}splf}(R)$.
\medskip

Moreover, we prove that the invariant $\mathrm{splf}(R)$ can be computed by determining the Gorenstein projective dimension, the projectively coresolved Gorenstein flat dimension, the Ding projective dimension, or the Gorenstein $\mathrm{AC}$-projective dimension of all flat left $R$-modules. This provides a refinement of \cite[Proposition 15]{DaleziosEmmanouil2023} due to Emmanouil and Dalezios. As a consequence, we obtain new characterizations of left $n$-perfect rings in terms of classes of Gorenstein modules (see Corollary~\ref{cor:Mejora_Prop3.2}). These characterizations extend and improve a result of Christensen, Estrada, and Thompson (see \cite[Proposition 3.2]{ChristensenEstradaThompson2021}), since we do not require the assumption that every flat Gorenstein projective left $R$-module is projective, a question that remains open (see \cite[Question 3.8]{BazzoniCortesIzurdiagaEstrada2020}). Furthermore, our results yield an extension of the second part of \cite[Theorem~4.9]{sarich2020singular} from the class of left perfect rings to the broader class of left $n$-perfect rings.

\medskip

\noindent \textbf{Theorem B.} Let $n \geq 0$ be an integer. For any ring $R$, the following statements are equivalent:
\begin{enumerate}
    \item Every flat left $R$-module has projective dimension at most $n$ (that is, $R$ is a left $n$-perfect ring).
    \item Every Gorenstein flat left $R$-module has Gorenstein projective dimension at most $n$.
    \item Every flat left $R$-module has Gorenstein projective dimension at most $n$.
   \item Every flat left $R$-module has Ding projective dimension at most $n$.
    \item Every Gorenstein flat left $R$-module has Ding projective dimension at most $n$.
    \item Every flat left $R$-module has $\mathrm{PGF}$ dimension at most $n$.
    \item Every Gorenstein flat left $R$-module has $\mathrm{PGF}$ dimension at most $n$.
    \item Every flat $R$-module has Gorenstein $\mathrm{AC}$-projective dimension at most $n$.
\end{enumerate}
Moreover, in each statement, the condition ``dimension at most $n$'' can be replaced by ``finite dimension''.
\medskip

Additionally, we introduce Gorenstein analogues of several classical homological invariants previously studied in \cite{EmmanouilTalelli2011}, namely $\mathrm{spli}(R)$, $\mathrm{silp}(R)$, $\mathrm{sfli}(R)$, and $\mathrm{silf}(R)$. We investigate the relationships between these invariants and their Gorenstein counterparts. In particular, we show that the Gorenstein invariants dominate their classical analogues (see Proposition~\ref{Prop:s<G-s}) and establish sufficient conditions under which they coincide (see Proposition \ref{prop:G-s+G-s=s+s}).

Finally, we characterize left weakly $n$-$\Sigma$-cotorsion rings, introduced in \cite{CortesIzurdiagaEstradaFresneda2026}, in terms of Gorenstein classes of modules, extending \cite[Theorem 5.1]{CortesIzurdiagaEstradaFresneda2026} (see Theorem~\ref{th:ch_of_weakly_n-Sigma_Cot_Rings}). Recall that a left $R$-module $M$ is \emph{$n$-$\Sigma$-cotorsion} if, for every set $I$, the direct sum $M^{(I)}$ has cotorsion dimension at most $n$ (for some integer $n \geq 0$). A ring $R$ is called \emph{left weakly $n$-$\Sigma$-cotorsion} if every injective left $R$-module is $n$-$\Sigma$-cotorsion. 

\medskip\noindent 
\textbf{Theorem C.}
Let $n \geq 0$ be an integer and $R$ any ring. The following statements are equivalent:
\begin{enumerate}
    \item Every injective left $R$-module is $n$-$\Sigma$-cotorsion (that is, $R$ is a left weakly $n$-$\Sigma$-cotorsion ring).
    \item Every Gorenstein injective left $R$-module is $n$-$\Sigma$-cotorsion.
    \item Every Ding injective left $R$-module is $n$-$\Sigma$-cotorsion.
    \item Every weakly Ding injective left $R$-module is $n$-$\Sigma$-cotorsion.
\end{enumerate}

The paper is organized as follows. In Section~\ref{section:Preliminaries}, we recall the necessary background on classical and Gorenstein homological algebra. Section~\ref{section:splf=G-splf} is devoted to the comparison between $\mathrm{splf}(R)$ and its Gorenstein analogue, where the main equality is established. In Section~\ref{section:Gorenstein invariants}, we introduce and study further Gorenstein invariants and investigate their relationships with the classical homological invariants.
Finally, in Section \ref{sect:weakly n-Sigma-Cotorsion rings}, we give characterizations of left weakly $n$-$\Sigma$-cotorsion rings introduced in \cite{CortesIzurdiagaEstradaFresneda2026} in terms of Gorenstein classes of modules.

\section{Preliminaries}\label{section:Preliminaries}
Throughout this paper, $R$ denotes an associative ring with identity, not necessarily commutative, and all $R$-modules are assumed to be unitary. We denote by $R^{op}$ the opposite ring of $R$. Unless explicitly stated otherwise, the term $R$-module always means a left $R$-module; the right $R$-modules will be specified explicitly when needed. We denote by $\Flat$, $\Inj$, $\Cot$, and $\Proj$ the classes of flat, injective, cotorsion, and projective $R$-modules, respectively. For an $R$-module $M$, we write $\fd_R(M)$, $\id_R(M)$, and $\pd_R(M)$ for the flat, injective, and projective dimensions of $M$, respectively.

Since some of our results involve definable classes, we briefly recall the relevant notion. A class of $R$-modules is called \emph{definable} if it is closed under direct products, direct limits, and pure submodules. Given a class $\mathcal{C}$ of $R$-modules, we denote by $\langle \mathcal{C} \rangle$ its \emph{definable closure}, i.e., the smallest definable class containing $\mathcal{C}$. In the case where $\mathcal{C} = \{C\}$, we simply write $\langle C \rangle$.

\subsection{Gorenstein and related classes of modules}

We recall the following standard definitions. An $R$-module $M$ is called:
\begin{enumerate}[label=(\arabic*)]
\item \emph{Gorenstein projective} (see \cite{EnochsJenda1995}) if $M$ occurs as a cycle in an exact complex of projective $R$-modules that remains exact under the functor $\mathrm{Hom}_R(-,P)$ for every projective $R$-module $P$. We denote by $\GP$ the class of all Gorenstein projective $R$-modules.
\item \emph{Gorenstein injective} (see \cite{EnochsJenda1995}) if $M$ occurs as a cycle in an exact complex of injective $R$-modules that remains exact under the functor $\mathrm{Hom}_R(E,-)$ for every injective $R$-module $E$. We denote by $\GI$ the class of all Gorenstein injective $R$-modules.
\item \emph{Gorenstein flat} (see \cite{EnochsJendaTorrecillas1993}) if $M$ occurs as a cycle in an exact complex of flat $R$-modules that remains exact after tensoring with any injective right $R$-module. We denote by $\GF$ the class of all Gorenstein flat $R$-modules.
\item \emph{Ding projective} (see \cite{ding2009strongly} and \cite{Gillespie2010}) if $M$ occurs as a cycle in an exact complex of projective $R$-modules that remains exact under the functor $\mathrm{Hom}_R(-,F)$ for every flat $R$-module $F$. We denote by $\DP$ the class of all Ding projective $R$-modules. It follows immediately from the definition that every Ding projective module is Gorenstein projective.
\item \emph{Ding injective} (see \cite{ding2009strongly} and \cite{Gillespie2010}) if $M$ occurs as a cycle in an exact complex of injective $R$-modules that remains exact under the functor $\mathrm{Hom}_R(T,-)$ for every $\FPInj$ $R$-module $T$. We denote by $\DI$ the class of all Ding injective $R$-modules. It follows immediately from the definition that every Ding injective module is Gorenstein injective.
\item \emph{weakly Ding injective} (see \cite{Iacob2023}) if $M$ occurs as a cycle in an exact complex of $\FPInj$ $R$-modules that remains exact under the functor $\mathrm{Hom}_R(T,-)$ for every $\FPInj$ $R$-module $T$. We denote by $\wDI$ the class of all weakly Ding injective $R$-modules. It follows immediately from the definition that every Ding injective module is weakly Ding injective.
\item \emph{projectively coresolved Gorenstein flat} (see \cite{sarich2020singular}) if $M$ occurs as a cycle in an exact complex of projective $R$-modules that remains exact after tensoring with any injective right $R$-module. We denote by $\mathrm{PGF}$ the class of all projectively coresolved Gorenstein flat $R$-modules. Moreover, by \cite[Corollary~4.5]{sarich2020singular}, one has $\PGF \subseteq \DP$ over any ring.
\item \emph{Gorenstein $\mathrm{AC}$-projective} (see \cite{BravoGillespieHovey2014}) if $M$ occurs as a cycle in an exact complex of projectives that remains exact under the functor $\mathrm{Hom}_R(-,L)$ for every level $R$-module $L$ (an $R$-module $L$ is \emph{level} if $\mathrm{Tor}_1^{R}(F,L)=0$ for each right $R$-module $F$ with a projective
resolution consisting of finitely generated projective right $R$-modules). We denote by $\GP_{ac}$ the class of all Gorenstein $\mathrm{AC}$-projective $R$-modules. Recall that, as a consequence of \cite[Proposition 2.10]{BravoGillespieHovey2014} together with \cite[Corollary 4.5]{sarich2020singular}, one has the inclusion $\GP_{ac} \subseteq \PGF$ over any ring.
\end{enumerate}

The notions of \emph{strongly Gorenstein projective}, \emph{strongly Gorenstein injective}, and \emph{strongly Gorenstein flat} modules were introduced by Bennis and Mahdou in \cite{stronglyBennisMahdou2007}. Recall that an $R$-module $M$ is called \emph{strongly Gorenstein projective} if it arises as a cycle of an exact complex of projective modules in which all components are equal and all differentials coincide, and such that this complex remains exact under the functor $\mathrm{Hom}_R(-,Q)$ for every projective module $Q$. The definitions of strongly Gorenstein injective and strongly Gorenstein flat modules are analogous. These classes are particularly important, as they provide direct-summand characterizations of the corresponding Gorenstein classes. More precisely, an $R$-module $M$ is Gorenstein projective if and only if it is a direct summand of a strongly Gorenstein projective module, and the same holds in the injective case (see \cite[Theorem 2.7]{stronglyBennisMahdou2007}). Moreover, as a consequence of \cite[Theorem~4.11]{sarich2020singular}, an $R$-module $M$ is Gorenstein flat if and only if it is a direct summand of a strongly Gorenstein flat module.

More generally, in \cite{BennisOuarghi2010}, Bennis and Ouarghi introduced a class of modules that generalizes the notion of Gorenstein projective modules. More precisely, let $\mathcal{X}$ be a class of $R$-modules that contains all projective modules. An $R$-module $M$ is said to be \emph{$\mathcal{X}$-Gorenstein projective} if $M$ occurs as a cycle in an exact complex of projective $R$-modules that remains exact under the functor $\mathrm{Hom}_R(-,X)$ for every module $X \in \mathcal{X}$. We denote by $\mathcal{X}\text{-}\GP$ the class of all $\mathcal{X}$-Gorenstein projective $R$-modules. This construction recovers several classes of modules defined previously, as summarized in the following remark.

\begin{remark} \label{remark: C-GP}
Let $\mathcal{X}$ be a class of $R$-modules that contains all projective modules.
    \begin{enumerate}
        \item If $\mathcal{X}=\Proj$, then $\mathcal{X}\text{-}\GP = \GP$ by definition.
        \item If $\mathcal{X}=\Flat$, then $\mathcal{X}\text{-}\GP = \DP$ by definition.
        \item If $\mathcal{X}=\langle R \rangle$, then $\mathcal{X}\text{-}\GP = \PGF$ by \cite[Corollary 4.5]{sarich2020singular}.
        \item If $\mathcal{X}$ is the class of all level $R$-modules, then $\mathcal{X}\text{-}\GP = \GP_{ac}$ by definition.
        \item If $\mathcal{X}$ is the class of all $R$-modules, then $\mathcal{X}\text{-}\GP = \Proj$ by definition.
    \end{enumerate}
\end{remark}

\subsection{Gorenstein and related homological dimensions}

Let $\mathcal{X}$ be a class of $R$-modules that contains all projective $R$-modules. The \emph{$\mathcal{X}$-Gorenstein projective dimension} of an $R$-module $M$, denoted by $\mathcal{X}\text{-}\Gpd_R(M)$, is defined as the least integer $n \geq 0$ such that there exists an exact sequence of $R$-modules
\[
0 \to G_n \to G_{n-1} \to \cdots \to G_0 \to M \to 0
\]
with each $G_i \in \mathcal{X}\text{-}\GP$. If no such integer exists, one sets $\mathcal{X}\text{-}\Gpd_R(M)= \infty$.

Observe that, by Remark \ref{remark: C-GP}, taking $\mathcal{X}$ to be the classes $\Proj$, $\Flat$, and $\langle R \rangle$ yields, respectively, the classical notions of \emph{Gorenstein projective dimension}, \emph{Ding projective dimension}, and \emph{projectively coresolved Gorenstein flat dimension} of $M$. These are denoted by $\Gpd_R(M)$, $\DPdim_R(M)$, and $\PGFdim_R(M)$, respectively. Likewise, the \emph{Gorenstein flat dimension} of $M$, denoted by $\Gfd_R(M)$, is defined analogously.

Dually, the \emph{Gorenstein injective dimension} of $M$, denoted by $\Gid_R(M)$, is defined as the least integer $n \geq 0$ such that there exists an exact sequence of $R$-modules
\[
0 \to M \to G^0 \to G^1 \to \cdots \to G^n \to 0
\]
with each $G^i \in \GI$. If no such integer exists, one sets $\Gid_R(M)=\infty$.

Using the previous notions, one can define Gorenstein analogues of the classical weak and global dimensions of a ring as follows:
\begin{enumerate}
    \item The \emph{left Gorenstein global dimension} of a ring $R$, denoted by $\mathrm{Ggldim}(R)$, is defined as the supremum of the Gorenstein projective dimensions (equivalently, by \cite[Theorem 1.1]{BennisMahdou2010}, the Gorenstein injective dimensions) of all $R$-modules. The \emph{right Gorenstein global dimension} of $R$ is defined as $\mathrm{Ggldim}(R^{op})$. 
    \item Similarly, the \emph{Gorenstein weak global dimension} of a ring $R$, denoted by $\mathrm{Gwgldim}(R)$, is defined as the supremum of the Gorenstein flat dimensions of all $R$-modules. Equivalently, by \cite[Corollary 2.5]{ChristensenEstradaThompson2021}, it can be computed as the supremum of the Gorenstein flat dimensions of all right $R$-modules.
\end{enumerate}

For any ring \(R\), there holds the inequality $\mathrm{Gwgldim}(R) \leq \mathrm{Ggldim}(R)$, as proved in \cite[Theorem~3.7]{WangYangShaoZhang2023}.

\begin{remark}\label{remark:resolutions}
Let $\mathcal{X}$ be a class of $R$-modules such that $\Proj \subseteq \mathcal{X}$. As shown in \cite{sarich2020singular} and \cite{BennisOuarghi2010}, the classes $\mathcal{X}\text{-}\GP$ and $\GF$ both contain all projective $R$-modules and are closed under kernels of epimorphisms, extensions, direct sums, and direct summands. Consequently, by \cite[Lemma 3.12]{auslander1969stable}, an $R$-module $M$ has $\mathcal{X}$-Gorenstein projective (respectively, Gorenstein flat) dimension at most $n$, for some integer $n \geq 0$, if and only if, for any partial projective resolution of $M$,
\[
0 \to G_n \to P_{n-1} \to \cdots \to P_1 \to P_0 \to M \to 0,
\]
the module $G_n$ lies in $\mathcal{X}\text{-}\GP$ (respectively, in $\GF$).

\noindent Dually, as shown in \cite{Holm2004Gorenstein}, the class $\GI$ contains all injective $R$-modules and is closed under cokernels of monomorphisms, extensions, direct products, and direct summands. It follows that an $R$-module $M$ has Gorenstein injective dimension at most $n$ if and only if, for any partial injective resolution of $M$,
\[
0 \to M \to E^0 \to E^1 \to \cdots \to E^{n-1} \to T^n \to 0,
\]
one has $T^n \in \GI$.
\end{remark}
\section{Projective vs. Gorenstein projective dimensions of \texorpdfstring{$\Flat$}{Flat} and \texorpdfstring{$\GF$}{Gorenstein Flat} modules} \label{section:splf=G-splf}
Let $R$ be a ring. We consider the following classical homological invariant:
\[
\mathrm{splf}(R) := \sup\{\pd_R(F) \mid F \in \Flat\}.
\]
The notation ``splf'' stands for \emph{supremum of projective lengths of flat modules}. Recall that $R$ is said to be a \emph{left $n$-perfect ring} if $\mathrm{splf}(R) \leq n$ for some integer $n\geq 0$. We also define the Gorenstein counterpart of the invariant $\mathrm{splf}(R)$ as follows:
\[
\mathrm{G\text{-}splf}(R) := \sup\{\Gpd_R(M) \mid M \in \GF\}.
\]
The following lemma is a slight improvement of \cite[Lemma~3.1]{ChristensenEstradaThompson2021}.
\begin{lemma}\label{lem:desigualdad_C-Gpd<Gfd+splf}
    Let $R$ be a ring and let $\mathcal{X}$ be a class of $R$-modules satisfying $\Proj \subseteq \mathcal{X} \subseteq \langle R \rangle$. Then, for any $R$-module $M$, the following inequality holds:$$\mathcal{X}\text{-}\Gpd_{R}(M) \leq \mathrm{Gfd}_{R}(M)+\sup \{ \mathcal{X}\text{-}\Gpd_{R}(L) \mid L \in \Flat \}.$$
\begin{proof}
First, observe that by Remark~\ref{remark: C-GP} one has $\PGF \subseteq \mathcal{X}\text{-}\GP \subseteq \GP$ in this case. Let $M$ be an $R$-module. The inequality is immediate if either $\Gfd_R(M)=\infty$ or $\sup \{ \mathcal{X}\text{-}\Gpd_{R}(L) \mid L \in \Flat \}=\infty$. Thus, we may assume that $\Gfd_R(M)=d<\infty$ and $\sup \{ \mathcal{X}\text{-}\Gpd_{R}(L) \mid L \in \Flat \}=n<\infty$. Under these assumptions, we prove that $\mathcal{X}\text{-}\Gpd_{R}(M) \leq d+n$.

\noindent Since $\Gfd_R(M)=d$, there exists an exact sequence:
\[
0 \to G_d \to L_{d-1} \to \cdots \to L_0 \to M \to 0,
\]
where each $L_i \in \Proj$, and $G_d\in \GF$ by Remark \ref{remark:resolutions}. Therefore, it is enough to show that $\mathcal{X}\text{-}\Gpd_{R}(G_d) \leq n$.

\noindent Since $G_d$ is Gorenstein flat, by \cite[Theorem 4.11]{sarich2020singular} there exists a short exact sequence: 
\[
0 \to G_d \to F \to N \to 0,
\] 
where $F\in \Flat$ and $N\in \PGF \subseteq \mathcal{X}\text{-}\GP$.

\noindent If \(n=0\), then $F \in \mathcal{X}\text{-}\GP$. Consequently, the preceding short exact sequence yields $G_d \in \mathcal{X}\text{-}\GP$, since the class $\mathcal{X}\text{-}\GP$ is closed under kernels of epimorphisms, as a consequence of \cite[Theorem 2.3(1)]{BennisOuarghi2010}.

\noindent Assume now that \(n>0\). By taking projective resolutions of \(G_d\) and \(N\), and applying the Horseshoe Lemma
(see \cite[Proposition~6.24]{Rotman2009}),
we obtain a commutative diagram with exact rows and columns in which all modules \(P_i\) and \(Q_j\) are projective \(R\)-modules:
\[
\vcenter{
\xymatrix{
& 0 \ar[d] & 0 \ar[d] & 0 \ar[d] & \\
0 \ar[r] & T_n \ar[r] \ar[d] & H_n \ar[r] \ar[d] & V_n \ar[r] \ar[d] & 0 \\
0 \ar[r] & P_{n-1} \ar[r] \ar[d] & P_{n-1} \oplus Q_{n-1} \ar[r] \ar[d] & Q_{n-1} \ar[r] \ar[d] & 0 \\
& \vdots \ar[d] & \vdots \ar[d] & \vdots \ar[d] & \\
0 \ar[r] & P_0 \ar[r] \ar[d] & P_0 \oplus Q_0 \ar[r] \ar[d] & Q_0 \ar[r] \ar[d] & 0 \\
0 \ar[r] & G_d \ar[r] \ar[d] & F \ar[r] \ar[d] & N \ar[r] \ar[d] & 0 \\
& 0 & 0 & 0 &
}
}
\]
\noindent Since $\sup \{ \mathcal{X}\text{-}\Gpd_{R}(L) \mid L \in \Flat \}=n$ and $F$ is flat, Remark~\ref{remark:resolutions} applied to the second column of the diagram implies that $H_n \in \mathcal{X}\text{-}\GP$. Furthermore, as $N \in \mathcal{X}\text{-}\GP$, the third column shows that $V_n \in \mathcal{X}\text{-}\GP$. It then follows from the top row of the diagram that $T_n \in \mathcal{X}\text{-}\GP$. Finally, the first column yields $\mathcal{X}\text{-}\Gpd(G_d) \leq n$.
\end{proof}
\end{lemma}

\begin{corollary}\label{cor:part_lem}
Let $R$ be a ring. Then, for any $R$-module $M$, one has
\[
\PGFdim_R(M) \leq \Gfd_R(M) + \sup \{ \PGFdim_R(L) \mid L \in \Flat \} \leq \Gfd_R(M) + \mathrm{splf}(R).
\]
\begin{proof}
The first inequality is a particular case of Lemma~\ref{lem:desigualdad_C-Gpd<Gfd+splf} obtained by taking $\mathcal{X}=\langle R \rangle$, in which case $\PGF = \mathcal{X}\text{-}\GP$ by Remark~\ref{remark: C-GP}. The second inequality follows from the fact that the projectively coresolved Gorenstein flat dimension is a refinement of the projective dimension.
\end{proof}
\end{corollary}

\begin{remark}
    The inequality in Lemma \ref{lem:desigualdad_C-Gpd<Gfd+splf} is, in general, strict. Indeed, it suffices to show that the first inequality in Corollary \ref{cor:part_lem}, which is a particular case of Lemma \ref{lem:desigualdad_C-Gpd<Gfd+splf}, may be strict. Consider the ring $\mathbb{Z}$ of integers. Since the global dimension of $\mathbb{Z}$ is equal to one, we have $\Proj=\PGF\subseteq\GF=\Flat$. Moreover, $\mathrm{splf}(\mathbb{Z})=1$. Hence, the first inequality in Corollary \ref{cor:part_lem} becomes $\pd_{\mathbb{Z}}(M)\leq \fd_{\mathbb{Z}}(M)+1$, which is strict for any non-flat $\mathbb{Z}$-module $M$.
\end{remark}

The following theorem shows that, over any ring $R$, one has $\mathrm{splf}(R)=\mathrm{G\text{-}splf}(R)$. Moreover, it recovers and extends \cite[Proposition 15]{DaleziosEmmanouil2023}, due to Emmanouil and Dalezios.

\begin{theorem}\label{th:s(GP)l(GF)=splf}
Let $R$ be a ring, and let $\mathcal{X}$ and $\mathcal{C}$ be classes of $R$-modules such that $\Proj \subseteq \mathcal{X} \subseteq \langle R \rangle$ and $\Proj \subseteq \mathcal{C}$. Then the following quantities coincide:
\begin{flalign*} 
s_1 &:= \sup \{ \mathcal{C}\text{-}\Gpd_R(M) \mid M \in \Flat \}, && \\ 
s_2 &:= \mathrm{G\text{-}splf}(R), && \\ 
s_3 &:= \mathrm{splf}(R), && \\ 
s_4 &:= \sup \{ \PGFdim_R(M) \mid M \in \GF \}, && \\
s_5 &:= \sup \{ \mathcal{X}\text{-}\Gpd_R(M) \mid M \in \GF \}.
\end{flalign*}

\begin{proof}
The proof proceeds by showing that: 
\[
s_2 \leq s_4 \leq s_3 \leq s_1 \leq s_3 \leq s_2,
\]
which immediately yields $s_1=s_2=s_3=s_4$. Once the preceding equalities have been established, the equality involving $s_5$ follows immediately from its definition.

\noindent To begin with, the inequality $s_2 \leq s_4$ is immediate from the inclusion $\PGF \subseteq \GP$, valid over any ring. Next, Corollary \ref{cor:part_lem} implies $s_4 \leq s_3$.

\noindent We now prove that $s_3 \leq s_1$. Assume that $s_1 = n < \infty$. By \cite[Theorem 4.2]{CortesIzurdiagaEstradaFresneda2026}, it suffices to show that for every set $I$, the module $R^{(I)}$ has cotorsion dimension at most $n$. Equivalently, by \cite[Proposition 19.2.1]{MaoDing2006CotorsionDimension}, one must verify that $\mathrm{Ext}_R^{n+1}(F, R^{(I)}) = 0$ for every flat $R$-module $F$.

\noindent If $n=0$, then every flat $R$-module is $\mathcal{C}$-Gorenstein projective. Since $R^{(I)}$ is projective and $\Proj \subseteq \mathcal{C}$, it follows immediately that $\mathrm{Ext}_R^{1}(F, R^{(I)})=0$ for all set $I$ and all flat $R$-module $F$.

\noindent Assume now that $n>0$, and let $F$ be a flat $R$-module. By Remark~\ref{remark:resolutions}, there exists an exact sequence
\[
0 \to K \to P_{n-1} \to \cdots \to P_0 \to F \to 0,
\]
where each $P_i \in \Proj$ and $K \in \mathcal{C}\text{-}\GP$. Dimension shifting yields $\mathrm{Ext}_R^{n+1}(F, R^{(I)}) \cong \mathrm{Ext}_R^{1}(K, R^{(I)})$. Since $K$ is $\mathcal{C}$-Gorenstein projective and $R^{(I)} \in \Proj$, the latter group vanishes. Hence $\mathrm{Ext}_R^{n+1}(F, R^{(I)})=0$, and therefore $s_3 \leq s_1$.

\noindent The reverse inequality $s_1 \leq s_3$ follows directly from the inclusion $\Proj \subseteq \mathcal{C}\text{-}\GP$. Moreover, the same argument used above to prove $s_3 \leq s_1$ 
also yields $s_3 \leq s_2$ in the special case $\mathcal{C}=\Proj$, since then $\mathcal{C}\text{-}\GP=\GP$ by Remark~\ref{remark: C-GP}. Consequently, $s_1=s_2=s_3=s_4$. Finally, since $\Proj \subseteq \mathcal{X} \subseteq \langle R \rangle$, then for any $R$-module $M$, one has $\Gpd_R(M) \leq \mathcal{X}\text{-}\Gpd_R(M) \leq \PGFdim_R(M)$. Consequently $s_2 \leq s_5 \leq s_4$, and therefore all five quantities coincide.
\end{proof}
\end{theorem}

\begin{remark}
The assumption in Theorem \ref{th:s(GP)l(GF)=splf} that 
$\Proj \subseteq \mathcal{X} \subseteq \langle R \rangle$ is required in order to guarantee the inclusion $\PGF \subseteq \mathcal{X}\text{-}\GP$. This inclusion is essential for the application of \cite[Theorem 4.11]{sarich2020singular} in the proof of Lemma \ref{lem:desigualdad_C-Gpd<Gfd+splf}. Since Corollary \ref{cor:part_lem} is a special case of this lemma and is subsequently used in the proof of Theorem \ref{th:s(GP)l(GF)=splf}, the above condition on $\mathcal{X}$ is essential throughout the argument.

\noindent We now present an example showing that if $\mathcal{X}$ is chosen to be the class of all $R$-modules (so that $\mathcal{X}\text{-}\GP=\Proj$ by Remark \ref{remark: C-GP}), then the corresponding supremum $\sup \{ \mathrm{pd}_R(M) \mid M \in \GF \}$ is strictly greater than $\mathrm{splf}(R)$. Let $p$ be a prime number and let $n \geq 3$ be an integer. Consider the ring $R = \mathbb{Z}/p^n\mathbb{Z}$, where $\mathbb{Z}$ denotes the ring of integers. By \cite[Corollary 3.10]{BennisMahdouOuarghi2010}, we have $\mathrm{Ggldim}(R)=0$, and consequently $\mathrm{Gwgldim}(R)=0$ by \cite[Theorem 3.7]{WangYangShaoZhang2023}. Since $R$ is not von Neumann regular while $\mathrm{Gwgldim}(R)=0$, it follows from \cite[Corollary 10.3.4]{EnochsJenda2011} that the weak global dimension of $R$ (and hence its global dimension) is infinite. Therefore, $\sup \{ \mathrm{pd}_R(M) \mid M \in \GF \} = \infty$. On the other hand, since $R$ is a perfect ring, we have that $\mathrm{splf}(R)=0$.
\end{remark}

\begin{corollary}\label{cor:Mejora_Prop3.2}
Let $n\geq0$ be an integer. The following conditions are equivalent over any ring $R$.
\begin{enumerate}
    \item Every flat left $R$-module has projective dimension at most $n$ (that is, $R$ is a left $n$-perfect ring).
    \item Every Gorenstein flat left $R$-module has Gorenstein projective dimension at most $n$.
    \item Every flat left $R$-module has Gorenstein projective dimension at most $n$.
   \item Every flat left $R$-module has Ding projective dimension at most $n$.
    \item Every Gorenstein flat left $R$-module has Ding projective dimension at most $n$.
    \item Every flat left $R$-module has $\mathrm{PGF}$ dimension at most $n$.
    \item Every Gorenstein flat left $R$-module has $\mathrm{PGF}$ dimension at most $n$.
    \item Every flat $R$-module has Gorenstein $\mathrm{AC}$-projective dimension at most $n$.
\end{enumerate}
Moreover, in each statement, the condition ``dimension at most $n$'' can be replaced by ``finite dimension''.
\begin{proof}
Applying Theorem~\ref{th:s(GP)l(GF)=splf} with $\mathcal{X}=\Flat$ and with $\mathcal{C}$ successively equal to $\Proj$, $\Flat$, $\langle R\rangle$, and the class of level $R$-modules. In each case, the corresponding supremum coincides with $\mathrm{splf}(R)$. Consequently, each of the above conditions is equivalent to the inequality $\mathrm{splf}(R)\leq n$, that is, to $R$ being a left $n$-perfect ring. The final assertion then follows from Remark~\ref{remark:resolutions} together with \cite[Lemma~3.9]{cortes2016products}, which shows that having finite dimension with respect to any of the classes considered is equivalent to the existence of a uniform finite bound.
\end{proof}
\end{corollary}

\begin{remark}
Corollary \ref{cor:Mejora_Prop3.2} improves and extends \cite[Proposition 3.2]{ChristensenEstradaThompson2021}, as it does not require the assumption that every flat and Gorenstein projective module is projective, a condition whose validity remains open, see \cite[Question 3.8]{BazzoniCortesIzurdiagaEstrada2020}. Moreover, Corollary \ref{cor:Mejora_Prop3.2} implies, in particular, that if $R$ is a ring with $\mathrm{splf}(R)=\infty$ (see \cite[Example 5.12]{CortesIzurdiagaEstradaFresneda2026} for examples of such rings), then there exists a flat $R$-module $F$ such that $\Gpd_R(F)=\infty$.
\end{remark}

Rings satisfying the conditions of Corollary~\ref{cor:Mejora_Prop3.2} for $n=0$ are precisely the left perfect rings introduced by Bass in \cite{Bass1960finitistic}. In this case, our results recover and extend the second part of \cite[Theorem~4.9]{sarich2020singular}, which characterizes left perfect rings as those for which the classes $\GF$ and $\PGF$ coincide. We record this special case separately in the following corollary. 

\begin{corollary}
The following conditions are equivalent over any ring $R$.
\begin{enumerate}
    \item $\Flat = \Proj$, that is, $R$ is a left perfect ring.
    \item $\GF = \PGF$.
    \item $\GF \subseteq \GP$.
    \item $\Flat \subseteq \GP$.
    \item $\GF \subseteq \DP$.
\end{enumerate}
\end{corollary}

\section{Gorenstein counterparts of classical homological invariants} \label{section:Gorenstein invariants}
In view of the results obtained in the previous section (Theorem~\ref{th:s(GP)l(GF)=splf}) concerning the relationship between the classical invariant $\mathrm{splf}(R)$ and its Gorenstein counterpart $\mathrm{G\text{-}splf}(R)$, we introduce in this section the Gorenstein analogues of other classical homological invariants and study some of their relationships. We first recall the standard notation for the classical invariants:\begin{align*}
\mathrm{spli}(R) &:= \sup \{ \pd_R(I) \mid I \in \Inj \}, \\
\mathrm{silp}(R) &:= \sup \{ \id_R(P) \mid P \in \Proj \}, \\
\mathrm{sfli}(R) &:= \sup \{ \fd_R(I) \mid I \in \Inj \}, \\
\mathrm{silf}(R) &:= \sup \{ \id_R(F) \mid F \in \Flat \}.
\end{align*}

The notation is motivated by the corresponding supremum of homological lengths, namely projective lengths of injective modules, injective lengths of projective modules, flat lengths of injective modules, and injective lengths of flat modules, respectively. We now introduce their Gorenstein analogues:
\begin{align*}
\mathrm{G\text{-}spli}(R) &:= \sup \{ \Gpd_R(E) \mid E \in \GI \}, \\
\mathrm{G\text{-}silp}(R) &:= \sup \{ \Gid_R(G) \mid G \in \GP \}, \\
\mathrm{G\text{-}sfli}(R) &:= \sup \{ \Gfd_R(E) \mid E \in \GI \}, \\
\mathrm{G\text{-}silf}(R) &:= \sup \{ \Gid_R(T) \mid T \in \GF \}.
\end{align*}
\begin{lemma}\label{lemma:Holm_lars}
Let $R$ be a ring and let $M$ be an $R$-module.
\begin{enumerate}
    \item If $\pd_R(M) < \infty$, then $\Gid_R(M) = \id_R(M)$.
    \item If $\id_R(M) < \infty$, then $\Gpd_R(M) = \pd_R(M)$ and $\Gfd_R(M) = \fd_R(M)$.
\end{enumerate}
\begin{proof}
Assertion (1) is proved in \cite[Theorem 2.1]{Holm2004}. 
Assertion (2) follows from \cite[Theorem 2.2]{Holm2004} and \cite[Corollary 2.2]{ChristensenEstradaThompson2021}.
\end{proof}
\end{lemma}

\begin{remark}\label{remark:desigualdades}
Let $R$ be a ring. Recall that the invariant
\[
\mathrm{FPD}(R):=\sup\{\operatorname{pd}_R(M) \mid M \text{ has finite projective dimension}\}
\]
is known as the \textit{finitistic projective dimension} of $R$. By a result of Jensen \cite[Proposition~6]{Jensen1970}, one has $\mathrm{splf}(R)\leq \mathrm{FPD}(R)$. Moreover, it follows from \cite[Proposition~1.3(i), Chapter~VII]{BeligiannisReiten2007} that $\mathrm{FPD}(R)\leq \mathrm{silp}(R)$, while \cite[Proposition~2.1]{EmmanouilTalelli2011} yields the equality $\mathrm{silp}(R)=\mathrm{silf}(R)$. Consequently, we have the following chain of inequalities:
\[
\mathrm{splf}(R)\leq \mathrm{FPD}(R)\leq \mathrm{silp}(R)=\mathrm{silf}(R).
\]
\end{remark}

The following proposition strengthens the well-known fact that the Gorenstein global dimension dominates the classical invariants: $\mathrm{spli}(R)$, $\mathrm{silp}(R)$, $\mathrm{silf}(R)$, and $\mathrm{sfli}(R)$.

\begin{proposition}\label{Prop:s<G-s}
Let $R$ be any ring. Then the following inequalities hold:
\begin{enumerate}
    \item $\mathrm{spli}(R)\leq\mathrm{G\text{-}spli}(R)$.
    \item $\mathrm{silp}(R)\leq \mathrm{G\text{-}silp}(R)$.
    \item $\mathrm{sfli}(R)\leq\mathrm{G\text{-}sfli}(R)$.
    \item $\mathrm{silf}(R)\leq\mathrm{G\text{-}silf}(R)$.
\end{enumerate}
\begin{proof}
We begin by proving (1). Assume that $\mathrm{G\text{-}spli}(R)=n<\infty$. Let $I \in \Inj$. By assumption, $\Gpd_R(I)\leq n$, and by Lemma \ref{lemma:Holm_lars} we have $\pd_R(I)=\Gpd_R(I)$. Hence, $\pd_R(I)\leq n$, and therefore $\mathrm{spli}(R)\leq n$.

\noindent We now prove (2). Assume that $\mathrm{G\text{-}silp}(R)=n<\infty$. Let $P \in \Proj$. By assumption, $\Gid_R(P)\leq n$, and by Lemma \ref{lemma:Holm_lars} we have $\id_R(P)=\Gid_R(P)$. It follows that $\id_R(P)\leq n$, and thus $\mathrm{silp}(R)\leq n$.

\noindent To prove (3), assume that $\mathrm{G\text{-}sfli}(R)=n<\infty$. Let $I \in \Inj$. By assumption, $\Gfd_R(I)\leq n$, and by Lemma \ref{lemma:Holm_lars} we have $\fd_R(I)=\Gfd_R(I)$. Consequently, $\fd_R(I)\leq n$, and hence $\mathrm{sfli}(R)\leq n$.

\noindent Finally, we prove (4). Assume that $\mathrm{G\text{-}silf}(R)=n<\infty$. In order to show that $\mathrm{silf}(R)\leq n$, it suffices to prove that $\mathrm{silp}(R)\leq n$, since these two quantities coincide by Remark \ref{remark:desigualdades}. Let $P \in \Proj$. By assumption, $\Gid_R(P)\leq n$, and by Lemma \ref{lemma:Holm_lars} we have $\id_R(P)=\Gid_R(P)$. Consequently, $\id_R(P)\leq n$, which implies that $\mathrm{silp}(R)\leq n$.
\end{proof}
\end{proposition}

The next proposition may be regarded as the Gorenstein analogue of \cite[Proposition 2.1]{EmmanouilTalelli2011}.

\begin{proposition} \label{prop:G-silf<G-silp}
For any ring $R$, one has $\mathrm{G\text{-}silf}(R) \leq \mathrm{G\text{-}silp}(R)$, with equality if every Gorenstein projective $R$-module is Gorenstein flat.
\begin{proof}
Assume that $\mathrm{G\text{-}silp}(R) = n$ for some integer $n\geq0$. By Theorem~\ref{th:s(GP)l(GF)=splf}, together with Remark~\ref{remark:desigualdades} and Proposition~\ref{Prop:s<G-s}, we obtain the chain of inequalities:
\[
\mathrm{G\text{-}splf}(R) = \mathrm{splf}(R) \leq \mathrm{silp}(R) \leq \mathrm{G\text{-}silp}(R) = n.
\]
\noindent Let $T \in \GF$. Then $\Gpd_R(T) \leq n$, and hence there exists an exact sequence:
\[
0 \to G_n \to G_{n-1} \to \cdots \to G_0 \to T \to 0,
\]
where $G_j \in \GP$ for all $j = 0, \dots, n$. Since $\mathrm{G\text{-}silp}(R)=n$, we have that $\Gid_R(G_j)\leq n$ for all $j=0, \dots, n$. Using that the class $\GI$ is closed under cokernels of monomorphisms (see \cite[Theorem 2.6]{Holm2004Gorenstein}), it follows that, for any integer $m \geq 0$, the class of modules with Gorenstein injective dimension at most $m$ is also closed under cokernels of monomorphisms. Hence, we deduce inductively that $\Gid_R(T) \leq n$. Consequently, $\mathrm{G\text{-}silf}(R) \leq n$. Finally, if every Gorenstein projective module is Gorenstein flat, the reverse inequality follows directly from the definitions.
\end{proof}
\end{proposition}

\begin{corollary}
Let $n \geq 0$ be an integer. The following statements hold for any ring $R$.

\begin{enumerate}
    \item If every Gorenstein injective $R$-module has Gorenstein projective dimension at most $n$, then every injective $R$-module has projective dimension at most $n$.

    \item If every Gorenstein projective $R$-module has Gorenstein injective dimension at most $n$, then:
    \begin{enumerate}
        \item every flat $R$-module has injective dimension at most $n$; and
        
        \item every Gorenstein flat $R$-module has Gorenstein injective dimension at most $n$.
    \end{enumerate}

    \item If every Gorenstein injective $R$-module has Gorenstein flat dimension at most $n$, then every injective $R$-module has flat dimension at most $n$.

    \item If every Gorenstein flat $R$-module has Gorenstein injective dimension at most $n$, then every flat $R$-module has injective dimension at most $n$.
\end{enumerate}
Moreover, in each statement, the condition ``dimension at most $n$'' can be replaced by ``finite dimension''.
\begin{proof}
The assertions follow immediately from Propositions~\ref{Prop:s<G-s} and~\ref{prop:G-silf<G-silp}. Moreover, the final statement follows from \cite[Lemma~3.9]{cortes2016products} and its dual in the injective setting, which show that finite dimension in any of the classes considered is equivalent to a uniform finite bound.
\end{proof}
\end{corollary}

The following theorem establishes a relationship between the classical invariants $\mathrm{silp}(R)$ and $\mathrm{spli}(R)$ and their Gorenstein counterparts $\mathrm{G\text{-}silp}(R)$ and $\mathrm{G\text{-}spli}(R)$. In particular, it provides a Gorenstein analogue of \cite[Proposition~1.3(vi), Chapter~VII]{BeligiannisReiten2007}: if $\mathrm{G\text{-}silp}(R) < \infty$ and $\mathrm{G\text{-}spli}(R) < \infty$, then these two invariants coincide.

\begin{proposition} \label{prop:G-s+G-s=s+s}
For any ring $R$, the following equalities hold:
\begin{enumerate}
    \item $\mathrm{silp}(R)+\mathrm{spli}(R)=\mathrm{G\text{-}silp}(R)+\mathrm{G\text{-}spli}(R)$.
    \item $\mathrm{silf}(R)+\mathrm{sfli}(R)=\mathrm{G\text{-}silf}(R)+\mathrm{G\text{-}sfli}(R)$.
\end{enumerate}
Moreover, if $\mathrm{silf}(R)+\mathrm{sfli}(R)<\infty$, then
\[
\mathrm{G\text{-}silp}(R)
=
\mathrm{silp}(R)
=
\mathrm{spli}(R)
=
\mathrm{G\text{-}spli}(R)
=
\mathrm{silf}(R)
=
\mathrm{G\text{-}silf}(R)
=
\mathrm{Ggldim}(R),
\]
and
\[
\mathrm{G\text{-}sfli}(R)
=
\mathrm{sfli}(R)
=
\mathrm{Gwgldim}(R).
\]
\begin{proof}
We first prove~(1). By Proposition~\ref{Prop:s<G-s} we have $\mathrm{silp}(R)+\mathrm{spli}(R)\leq\mathrm{G\text{-}silp}(R)+\mathrm{G\text{-}spli}(R)$, so it suffices to establish the reverse inequality. Assume that $\mathrm{silp}(R)+\mathrm{spli}(R)<\infty$. Then, by \cite[Theorem 2.2, Chapter VII]{BeligiannisReiten2007}, one has $\mathrm{silp}(R)=\mathrm{spli}(R)=\mathrm{Ggldim}(R)<\infty$. Applying Proposition~\ref{Prop:s<G-s} again yields
\[
\mathrm{G\text{-}silp}(R)
=
\mathrm{silp}(R)
=
\mathrm{spli}(R)
=
\mathrm{G\text{-}spli}(R),
\]
which proves (1).

\noindent We now prove~(2). Proposition~\ref{Prop:s<G-s} gives $\mathrm{silf}(R)+\mathrm{sfli}(R)\leq\mathrm{G\text{-}silf}(R)+\mathrm{G\text{-}sfli}(R)$, so it remains to show the reverse inequality. Assume that $\mathrm{silf}(R)+\mathrm{sfli}(R)<\infty$. By \cite[Corollary 3.13]{WangYangShaoZhang2023}, it follows that $\mathrm{silp}(R)+\mathrm{spli}(R)<\infty$, and hence, by the argument established in the proof of~(1), we obtain
\[
\mathrm{G\text{-}silp}(R)
=
\mathrm{silp}(R)
=
\mathrm{spli}(R)
=
\mathrm{G\text{-}spli}(R)
=
\mathrm{Ggldim}(R).
\]
Combining these equalities with Proposition~\ref{Prop:s<G-s}, taking into account that $\mathrm{silp}(R)=\mathrm{silf}(R)$ over any ring by Remark \ref{remark:desigualdades}, we obtain:
$\mathrm{Ggldim}(R) = \mathrm{silp}(R) = \mathrm{silf}(R) \leq \mathrm{G\text{-}silf}(R) \leq \mathrm{Ggldim}(R),$ 
and thus all inequalities are equalities. Hence, $\mathrm{silf}(R)=\mathrm{G\text{-}silf}(R)=\mathrm{Ggldim}(R)$. Finally, \cite[Corollary 3.13]{WangYangShaoZhang2023} also implies that $\mathrm{sfli}(R^{\mathrm{op}})<\infty$, and therefore, by \cite[Theorem 5.3]{Emmanouil2012}, $\mathrm{sfli}(R)=\mathrm{sfli}(R^{\mathrm{op}})=\mathrm{Gwgldim}(R)$. Applying Proposition \ref{Prop:s<G-s} once more gives $\mathrm{Gwgldim}(R)=\mathrm{sfli}(R)\leq\mathrm{G\text{-}sfli}(R)\leq\mathrm{Gwgldim}(R),$ so all inequalities are equalities and $\mathrm{G\text{-}sfli}(R)=\mathrm{sfli}(R)=\mathrm{Gwgldim}(R)$. This completes the proof of (2) and the final assertion.
\end{proof}
\end{proposition}

\section{Gorenstein characterizations of weakly \texorpdfstring{$n\text{-}\Sigma$}{n-Sigma}-cotorsion rings}\label{sect:weakly n-Sigma-Cotorsion rings}
In this section, we obtain characterizations of left weakly \(n\)-\(\Sigma\)-cotorsion rings in terms of Gorenstein classes of modules, extending \cite[Theorem 5.1]{CortesIzurdiagaEstradaFresneda2026}. 

Let $R$ be a ring. An \(R\)-module \(M\) is said to be \emph{\(n\)-\(\Sigma\)-cotorsion} provided that, for every set \(I\), the direct sum \(M^{(I)}\) has cotorsion dimension at most \(n\), where \(n \geq 0\) is an integer (see \cite[Definition 3.1]{CortesIzurdiagaEstradaFresneda2026}). A ring \(R\) is called \emph{left weakly \(n\)-\(\Sigma\)-cotorsion} if every injective \(R\)-module is \(n\)-\(\Sigma\)-cotorsion (see \cite[Definition 5.2]{CortesIzurdiagaEstradaFresneda2026}). The class of left weakly $n$-$\Sigma$-cotorsion rings includes both the class of left noetherian rings and the class of left $n$-perfect rings (see \cite[Example 5.10]{CortesIzurdiagaEstradaFresneda2026}). Moreover, it also contains rings that are neither left noetherian nor left perfect (see \cite[Example 5.16(2)]{CortesIzurdiagaEstradaFresneda2026}).

To this end, we recall the notion of periodicity with respect to a class of modules, which will be used throughout. Following \cite{BazzoniCortesIzurdiagaEstrada2020}, let \(\mathcal{C}\) be a class of \(R\)-modules. An \(R\)-module \(M\) is said to be \emph{\(\mathcal{C}\)-periodic} if there exists a short exact sequence $0 \to M \to C \to M \to 0$, with \(C \in \mathcal{C}\), and it is called \emph{trivial} if it belongs to $\mathcal{C}$.

For any integer $n\geq0$, we write $\Cot_n$ for the class of $R$-modules of cotorsion dimension at most $n$, noting that $\Cot_0=\Cot$. The following lemma extends the second part of \cite[Proposition 4.8(2)]{BazzoniCortesIzurdiagaEstrada2020} to the class $\Cot_n$, and will be essential in the proof of the main result of this section.

\begin{lemma}\label{lema: Every Cot_n-periodic is trivial}
Let $R$ be a ring and let $n \geq 0$ be an integer. Then every $\Cot_n$-periodic $R$-module belongs to $\Cot_n$, that is, every $\Cot_n$-periodic $R$-module is trivial.
\begin{proof}
The case $n=0$ is precisely \cite[Proposition 4.8(2)]{BazzoniCortesIzurdiagaEstrada2020}. Assume therefore that $n>0$, and let $M$ be a $\Cot_n$-periodic $R$-module. Thus, there exists a short exact sequence $0 \to M \to T \to M \to 0$, with $T \in \Cot_n$. Consider a partial injective resolution of $M$. By applying the Horseshoe Lemma (see \cite[Proposition 6.24]{Rotman2009}) to the above short exact sequence, one obtains a commutative diagram with exact rows and columns:
\[
\vcenter{
\xymatrix{ & 0 \ar[d] & 0 \ar[d] & 0 \ar[d] & \\ 0 \ar[r] & M \ar[r] \ar[d] & T \ar[r] \ar[d] & M \ar[r] \ar[d] & 0 \\ 0 \ar[r] & E^{0} \ar[r] \ar[d] & E^{0} \oplus E^{0} \ar[r] \ar[d] & E^{0} \ar[r] \ar[d] & 0 \\ & \vdots \ar[d] & \vdots \ar[d] & \vdots \ar[d] & \\ 0 \ar[r] & E^{n-1} \ar[r] \ar[d] & E^{n-1} \oplus E^{n-1} \ar[r] \ar[d] & E^{n-1} \ar[r] \ar[d] & 0 \\ 0 \ar[r] & K^n \ar[r] \ar[d] & C^n \ar[r] \ar[d] & K^n \ar[r] \ar[d] & 0 \\ & 0 & 0 & 0 & }
}
\]
where each $E^i$ is an injective $R$-module.

\noindent Since $T \in \Cot_n$, the exactness of the second column implies that $C^n \in \Cot$. Consequently, the bottom row yields a short exact sequence
$0 \to K^n \to C^n \to K^n \to 0$ showing that $K^n$ is a $\Cot$-periodic $R$-module, then \cite[Proposition 4.8(2)]{BazzoniCortesIzurdiagaEstrada2020} implies that $K^n \in \Cot$. Finally, the exactness of the first column shows that $M \in \Cot_n$, as required.
\end{proof}
\end{lemma}

\begin{theorem}\label{th:ch_of_weakly_n-Sigma_Cot_Rings}
Let $n \geq 0$ be an integer and let $R$ be a ring. The following statements are equivalent:
\begin{enumerate}
    \item Every injective $R$-module is $n$-$\Sigma$-cotorsion, that is, $R$ is a left weakly $n$-$\Sigma$-cotorsion ring.
    
    \item Every Gorenstein injective $R$-module is $n$-$\Sigma$-cotorsion.
    
    \item Every Ding injective $R$-module is $n$-$\Sigma$-cotorsion.
    
    \item Every weakly Ding injective $R$-module is $n$-$\Sigma$-cotorsion.
    
    \item Every $\langle \Inj \rangle$-periodic $R$-module is $n$-$\Sigma$-cotorsion.
\end{enumerate}

\begin{proof}
By \cite[Theorem 2.7]{stronglyBennisMahdou2007}, every Gorenstein injective $R$-module is a direct summand of a strongly Gorenstein injective $R$-module. Since every strongly Gorenstein injective $R$-module is, in particular, $\langle \Inj \rangle$-periodic, it follows that every Gorenstein injective $R$-module is a direct summand of a $\langle \Inj \rangle$-periodic $R$-module. Moreover, the argument employed in the proof of \cite[Theorem 2.7]{stronglyBennisMahdou2007} also shows that every Ding injective and every weakly Ding injective $R$-module is a direct summand of a $\langle \Inj \rangle$-periodic $R$-module (for the class of weakly Ding injective modules, it suffices to observe that the class of $\mathrm{FP}$-injective $R$-modules is contained in $\langle \Inj \rangle$). Since the class $\Cot_n$ is closed under direct summands by \cite[Proposition 19.2.1]{MaoDing2006CotorsionDimension}, condition~(5) implies each of the remaining conditions. On the other hand, it is immediate that each of conditions~(2), (3), and~(4) implies condition~(1). Therefore, in order to prove the equivalence of the five conditions, it suffices to show that~(1) implies~(5).

\noindent Assume that condition~(1) holds, and let $M$ be a $\langle \Inj \rangle$-periodic $R$-module. Then there exists a short exact sequence $0 \to M \to E \to M \to 0$, where $E \in \langle \Inj \rangle$. Let $I$ be any set. Taking direct sums indexed by $I$, we obtain a short exact sequence $0 \to M^{(I)} \to E^{(I)} \to M^{(I)} \to 0$. Since definable classes are closed under arbitrary direct sums, it follows that $E^{(I)} \in \langle \Inj \rangle$. Now, by condition~(1) together with \cite[Theorem 5.1]{CortesIzurdiagaEstradaFresneda2026}, we obtain $\langle \Inj \rangle \subseteq \Cot_n$. Therefore, $E^{(I)} \in \Cot_n$, and consequently $M^{(I)}$ is a $\Cot_n$-periodic $R$-module. Applying Lemma~\ref{lema: Every Cot_n-periodic is trivial}, we conclude that $M^{(I)} \in \Cot_n$. Thus, $M$ is $n$-$\Sigma$-cotorsion.
\end{proof}
\end{theorem}

\section*{Acknowledgements}
The author would like to express his sincere gratitude to his PhD advisors, Sergio Estrada and Manuel Cortés-Izurdiaga, for their guidance, support, and many valuable comments and observations throughout the preparation of this manuscript.

\printbibliography[heading=bibintoc]
\end{document}